\documentstyle{amsart}

\title[Exactly k-to-1  maps]{Exactly k-to-1  maps: from pathological functions 
with
finitely many discontinuities to well-behaved covering maps}
\author{Jo Heath}
\address{Mathematics Department, Auburn University, Auburn, AL36849-5310}
\email{heathjw@@mail.auburn.edu}
\keywords{2-to-1 map, k-to-1 map, continua, graphs, dendrites
 indecomposable spaces, covering maps, simple maps, finitely discontinuous
 maps, tree-like continua}
\subjclass{primary 54C10, secondary 26A03}

\begin{document}

\setlength{\baselineskip}{14pt}

\maketitle

\vskip .5in

\def\picture #1 by #2 (#3){
  \vbox to #2{
    \hrule width #1 height 0pt depth 0pt
    \vfill
    \special{picture #3} 
    }
  }

\def\map {\picture 4.01in by 1.71in (map scaled 1000)}

\hrule
\vskip .1in
\noindent {\sc Abstract.} Many mathematicians encounter k-to-1  maps
only in the study of covering maps. But, of course, k-to-1 maps do
not have to be open. This paper touches on covering maps, and simple maps, but 
concentrates on ordinary k-to-1 functions (both continuous and
finitely discontinuous) from one metric continuum to another. New results, old 
results, and ideas for further research are given; and a baker's dozen of
questions are raised.
\vskip .1in
\hrule
\vskip .1in

\section{Introduction}
Requiring a function from one metric continuum $X$ to another, $Y$, to be 
finite-to-one, or even to be light, adds a strong hypothesis. But if the 
function must be k-to-1, meaning that each inverse has exactly $k$ points,
then the collection of available maps shrinks drastically and may even 
disappear. For instance,  if $Y$ is a dendrite, then there is a wealth of 
finite-to-one maps that map onto $Y$, but there are no  k-to-1 maps,
\cite{gotts}. What is it about the dendrite that repels these maps? Now, 
consider the domain $X$. It may be that every metric continuum $X$ admits a
k-to-1 map for $ k > 2 $ (see Question 9 later), but for the special
case $ k = 2 $ many interesting situations arise. For instance, the unit 
interval does not admit an exactly 2-to-1 map, \cite{harrold1}, but
some dendrites do. What is the crucial topological difference between an arc 
and a dendrite?  And one of the big questions today in this field is whether or 
not the pseudo-arc admits such a map. The central purpose of this survey paper 
is to describe what is known about the domains and images of exactly
k-to-1  maps, with special emphasis on the important $ k = 2 $ case, and to 
list
some of the many questions that still need to be answered. Secondly we will see 

what happens when finitely many discontinuities are allowed; the surprising 
thing is that there is still a lot of control. Thirdly and fourthly we will 
touch on two related topics: covering maps when the spaces do not have the 
usual textbook connectedness properties and, in variance with the title of this 

survey,  simple maps, that is, maps such that each point in the image has an 
inverse of cardinality 1 or 2. 

By {\em continuum} we mean a connected compact metric space. There is a 
glossary at the end of the paper containing other definitions. 
 
\section{Domains and images of k-to-1 maps.}

\subsection{2-to-1 images.}
If one toys with the question of which continua are 2-to-1 images
(of continua), it is quickly seen that it is easy to map 2-to-1 onto
a circle and onto other continua with subcontinua that are similar (in some 
sense) to a circle and it is hard otherwise. In fact, in \cite{nw}, Nadler and 
Ward show how to construct a straightforward 2-to-1 map (or k-to-
1  for any $k > 1$) onto any continuum that contains a non-unicoherent
subcontinuum. In this same paper, they ask if any tree-like continuum could be 
the image of an exactly 2-to-1 map. Their question, still unanswered
today, is the basis of the conjecture that a continuum is a 2-to-1
image iff it is not tree-like. Furthermore, each of the examples of non-tree-
like continua, that the author has tested, is a 2-to-1
retract of some continuum. Hence the following two questions, if answered 
affirmatively by a helpful reader, would neatly classify continua that
are  2-to-1 images.
\vspace{.2in}
{\bf Question 1.} \cite{nw} Is it true that no tree-like continuum 
can be the 2-to-1 image of a continuum?
\vspace{.2in}

{\bf Question 2.} Is it true that every continuum that is not tree-like is the 
2-to-1 image of a continuum? Can ``image"  be replaced with
``retract"?
\vspace{.2in}

Regarding Question 1, we know that the following types of continua, 
if tree-like,  cannot be the 2-to-1 image of any continuum:
dendrites (\cite{gotts}), hereditarily indecomposable continua 
(\cite{hiimage}), and indecomposable arc-continua (\cite{jdm1}). Furthermore, 
if a continuum has any of the following properties, then it 
cannot be the  2-to-1 image of a continuum: (1) every subcontinuum
has a cut point, \cite{nw} and \cite{gotts}, (2) every subcontinuum has a 
finite separating set and the continuum is hereditarily unicoherent, 
\cite{tree}, or (3) every subcontinuum has an endpoint, \cite{nw}.

Furthermore, we know that whatever the tree-like continuum $Y$, there is no 
confluent or crisp (\cite{hiimage}) 2-to-1 map from any continuum
onto $Y$.

Regarding Question 2, the 2-to-1 maps onto orientable or non-orientable
indecomposable arc-continua that are local Cantor bundles (which includes all 
solenoids for instance) have recently been studied in \cite{jdm2}. (These 
definitions are in the glossary.) It was proved that in the non-orientable 
case, every 2-to-1 map onto the continuum is a 2-fold covering map and in
either case every 2-to-1 map onto the continuum is either a 2-fold
cover or a retraction. Furthermore, every orientable local Cantor bundle is 
the  2-to-1 image of a continuum.
 
\subsection{k-to-1  images.}
If we consider integers larger than $ k = 2 $, the situation is murky. For each 

of these larger $ k$, there is indeed a tree-like continuum that is the
k-to-1 image of a continuum, \cite{tree}. On the other hand, no
dendrite is the k-to-1 image of a continuum for any $ k > 1 $,
\cite{gotts}. There isn't a lot of elbow room between dendrites and 
tree-like continua, and we do not even have a conjecture as to what the 
classification might be: 
\vspace{.2in}
{\bf Question 3.} For integers $ k > 2 $, which continua are k-to-1
images?
\vspace{.2in}

A related question asks which continua are k-to-1 images of
dendrites \cite{miklos1}. The topological structure of a dendrite dictates, 
\cite{dendrite}, that any k-to-1 image must be one-dimensional, it
must contain a simple closed curve, and it cannot contain uncountably many 
disjoint arcs. And of course it must be a Peano continuum. Is this sufficient? 
Yes, {\em if} the continuum contains only finitely many simple closed curves  ( 

\cite{dendrite}, \cite{miklos1}); but sometimes the answer is yes when the 
continuum does contain infinitely many simple closed curves \cite{dendrite}. 
\vspace{.2in}
{\bf Question 4.} \cite{miklos1}  Exactly which continua are k-to-1
images of dendrites?
\vspace{.2in}

{\bf Question 5.} \cite{dendrite} If each of $Y_{1}$ and $Y_{2}$ is the 
k-to-1 image of a dendrite, and if $Y_{1}  \cap Y_{2}$ is a single
point, then must $Y_{1}  \cup   Y_{2}$ be the k-to-1 image of a
dendrite?
\vspace{.2in}
{\bf Question 6.} \cite{dendrite} Might the answer to Question 4 depend on $k$?

That is, does there exist a continuum $Y$ and integers $k$ and $m$, both 
greater than 1, such that $Y$ is the k-to-1 image of a dendrite, but
$Y$ is not the m-to-1  image of a dendrite?
\vspace{.2in}
\subsection{2-to-1  domains}
It has been known for over fifty years that no 2-to-1 map can be
defined on an arc \cite{harrold1}, or, in fact,  any connected graph with odd 
Euler number  \cite{gilbert}, but extending these results has been difficult. 
Many dendrites admit 2-to-1 maps, and Wayne Lewis, \cite{tree}, has
constructed a decomposable arc-like continuum that admits a 2-to-1
map. In \cite{debski1} W. D\c{e}bski uses the  time-honored technique (see for 
instance \cite{mio}, \cite{cern1} and \cite{civin}) of classifying  continuous
involutions on a space in order to indirectly study 2-to-1 maps on the
space. He applies this to determine that a solenoid admits a 2-to-1
map iff there are at most finitely many even integers in an integer sequence 
that defines the inverse limit structure of the solenoid.

It would be nice of course to know exactly which continua admit 2-
to-1  maps but we don't even know the answers to these more restricted
questions:

\vspace{.2in}
{\bf Question 7} \cite{tree} Is there an indecomposable arc-like continuum that 

admits a
2-to-1 map?
\vspace{.2in}
(No use trying the classic Knaster Buckethandle space; J. Mioduszewski proved 
\cite{mio} over thirty years ago that it does not admit a 2-to-1
map. See also \cite{debski1}.)
\vspace{.2in}
{\bf Question 8} \cite{mio} Does the pseudo-arc admit a 2-to-1 map?
\vspace{.2in}
It is known  \cite{weak} that there is no weakly confluent 2-to-1
map defined on the pseudoarc. In fact, if $f$ is a weakly confluent 2-
to-1  map defined on any hereditarily indecomposable continuum, then neither
the domain or the image can be tree-like. For more information related to 
Question 8, see \cite{hidomain}.

George Henderson \cite{hend} has proved that if the domain, $ X $, of a
2-to-1 map is a mod 2 homology sphere  and $ \phi(X) $ denotes the
homology dimension of $ X $, then $ H_{i}(Y,Z_{2}) = Z_{2} $ if $ 0 \leq i \leq 
\phi(X) $ and $ H_{i}(Y,Z_{2}) =  0 $ otherwise, where $ Y $ is the image. Thus 
additive mod 2 homology cannot be used to distinguish 2-to-1 images of such a
sphere. As a corollary, he has that the circle is the only sphere that maps 2-
to-1  onto a sphere.
 
\subsection{k-to-1  domains}
As is evident by the question below, very little is known about which continua 
can be the domain of a k-to-1 map, if $ k > 2 $.
. 
\vspace{.2in}
{\bf Question 9.} Is there a continuum $X$ and an integer $k > 2 $
such that there is no exactly k-to-1 map defined on $X$?
\vspace{.2in}
There is extensive literature concerning which continua can or cannot be 
covering spaces, i.e. domains of very special k-to-1 maps, namely
covering maps. We will not attempt to survey these results, but we will mention 

two relatively recent papers. R. Myers \cite{myer} has constructed contractible 

open 3-manifolds which cannot cover closed 3-manifolds; and David Wright 
\cite{wright} gives a general method of determining when a contractible 
manifold cannot be a covering space of a manifold.
 
\subsection{k-to-1 maps between graphs.}
There should be some way to look at the adjacency matrices of two given graphs 
and decide if a k-to-1 map exists from one onto the other. Although
good progress has been made on this question, a direct answer has not been 
found (see Question 10 below). In this discussion we assume that the two given 
graphs have enough vertices to eliminate loops, are non-trivial (consist of 
more than just one vertex), and are connected, even though many of the known 
results are true, with little or no modification, for disconnected graphs. 

Given a positive integer $k$, and graphs $G$ and $H$, there are some 
preliminary filters to rule out the existence of a k-to-1 map from $G$
onto $H$.  For instance, is $k$ times the Euler number of $H$ at least as large 

as the Euler number of $G$ (or, if $ k = 2 $, is the Euler number of $G$ twice 
that of $H$)? If not, then there 
is no finitely discontinuous  k-to-1 function from $G$ onto $H$, much
less a continuous one. (See the theorem stated later in the subsection on
finitely discontinuous functions.) Since only endpoints (vertices of order one) 

of $G$ can map to endpoints of $H$, one can count them and make sure $G$ has at 

least $k$ times as many as $H$. A more subtle requirement, true for odd 
integers $k$, is that each vertex of $H$ with odd order must have an odd number 

of vertices in $G$ with odd order mapping to it (\cite{johilton1} or 
\cite{hilton2}). So, one can make sure that $G$ has at least as many odd-order 
vertices as $H$ does. But these tests can only give a definite ``no". S. Miklos 

has one of the few definite yes results in \cite{miklos2}; namely, if $k$ is 
odd, then a graph admits a k-to-1 map onto itself iff it has no
endpoints.

The original paper \cite{johilton1} that worked on Question 10 started with 
graphs $G$ and $H$ and a k-to-1 function $f$ from a vertex set of $G$
onto a vertex set of $H$, and answered whether or not $ f$ extended to a 
k-to-1 map from all of $G$ onto all of $H$. Similar questions for
$ \leq$k-to-1 maps are answered in \cite{johilton3} and
\cite{johilton4}.   The answers are algebraic in terms of the 
adjacency matrix for $H$ and the ``inverse adjacency" matrix for $G$ and $f$ 
(defined in the glossary).  An example of one of the theorems is as follows:
\vspace{.2in}
Theorem. \cite{johilton1} Suppose $G$ and $H$ are graphs, $k$ is an odd 
integer, and $f$ is a k-to-1 function from a vertex set of $G$ onto a
vertex set of $H$.  Then $f$ extends to a k-to-1 map from all of $G$
onto $H$ iff $f$, the adjacency matrix $A$, and the inverse adjacency matrix 
$B$ satisfy:
\begin{enumerate}
\item For each vertex $p$ in $H$, $k$ times the order of $p$ is at least as 
large as the sum of the orders of the vertices in $f^{-1}(p)$,
\item each diagonal element of $ k \cdot A - B $ is even and non-negative, and 
\item each entry of $ B - A $ is nonnegative.
\end{enumerate}

\vspace{.2in}
The shortcoming of this theorem and the other similar results is clear. If a 
given k-to-1 function from the vertex set of $G$ onto the vertex set
of $H$ fails to extend, that does {\em not} mean that there is no k-to-
1 function from $G$ onto $H$. Perhaps we just started with the wrong vertex
function.

A favorite approach is to change the question. Given two non-trivial, connected 

graphs, $G$ and $H$, does there {\em exist} an integer $k$ (or odd integer  $k$ 

or even integer $k$) and a k-to-1 map from $G$ onto $H$ \cite{hilton1}?
Another variation is: suppose that $G$ and $H$ are compatible enough to 
admit a k-to-1 map, and suppose $m$ is a larger integer (perhaps with
the same parity);  must $G$ and $H$ admit a m-to-1 map
\cite{johilton5} ? Many cases of these and other similar questions  have 
affirmative answers and the answers depend loosely on how close $H$ is to being 

a simple closed curve.  If $ H \neq S^{1}$, then $H$ is inspected to see if it 
is at least Eulerian (every vertex has even order). If not, how many odd-order 
vertices does $H$ have compared to the number of odd-order vertices in $G$; 
and, most important, how many endpoints does $H$ have? The graph $G$ seems to 
have little to do with the answer in many cases. Three nice results by A.J.W. 
Hilton \cite{hilton1} are: (1) if $H = S^{1}$ and $k$ is greater than the 
number of vertices in $G$, then all that is needed for the existence of a 
k-to-1 map from $G$ onto $H$ is that $G$ not have  more edges than
vertices. (2) If $H \neq S^{1}$ but $H$ has no endpoints, then for all 
sufficiently large even integers $k$, there is a k-to-1 map from $G$
onto $H$, and (3) if $H \neq S^{1}$, $H$ has no endpoints and there are at 
least as many odd-order vertices in $G$ as there are in $H$, then there are k-
to-1 maps from $G$ onto $H$ for all sufficiently large odd $k$. The
conditions given for $G$ and $H$ are, in each of the three cases, also 
necessary. Hilton \cite{hilton2} has also studied the relationships, for each 
parity, between  
the {\em initial} $k$ (the least $k$ such that there is a k-to-1 map) and the
{\em threshold} $k$ (the least $k$ such that $k$ and every integer larger than 
$k$
admits a k-to-1 map) and found that, for each parity, in many cases they are 
the same
integer. See \cite{hilton3} and \cite{johilton2} for some explicit 
constructions of k-to-1 maps from graphs onto a simple closed curve.

\vspace{.2in}
{\bf Question 10.} Given an integer $k > 1$ and two graphs, $G$ and $H$, when 
does a k-to-1 map exist from $G$ onto $H$?
\vspace{.2in}
\section{Finitely discontinuous k-to-1 functions.}

With many studies involving functions, so many theorems go out the window if a 
discontinuity is allowed that much of the power is lost. But not so with k-to-
1  functions, especially when the image is required to be a continuum. (Even
with k-to-1 functions, a single discontinuity can easily destroy
both connectivity and compactness.) The process is this: Suppose the domain is 
a graph $G$. Remove a finite set, $ N$, of points from $G$. Now reassemble the 
components of the complement of $N$, along with the points of  $ N$, in a
k-to-1 fashion in such a way that the resulting space is a continuum.
From this mental picture emerges the fact that the Euler number of $G$ is 
all-important. (We define the Euler number of a graph to be the number of edges 

minus the number of vertices.) In fact, the following theorem is a concise 
characterization of exactly which pairs of graphs have k-to-1
finitely discontinuous functions between them. No such characterization, based 
entirely on Euler numbers, is possible for continuous k-to-1
functions between graphs; in fact, we have no characterization at all for the 
continuous case (see earlier subsection). So, in the case of graphs, allowing
finitely many discontinuities actually clarifies the picture. For some studies 
of finitely discontinuous functions from an arc {\em into} an arc, where the 
image is not required to be compact, see \cite{katkel}, \cite{kat} and 
\cite{even}.
\vspace{.2in}
Theorem. \cite{fdgraph} If $G$ is a graph with Euler number $ m $ and $ H $ is 
a graph with Euler number $ n $, then there is a k-to-1 function
from $G$ onto $H$ with finitely many discontinuities:
\begin{enumerate}
\item iff $ m \leq kn $, if $ k > 2 $, and
\item iff $ m = 2n $, if $ k = 2 $.
\end{enumerate}
\vspace{.2in}
Up to now, the study of finitely discontinuous k-to-1 functions has
remained mostly in the safe haven of locally connected continua for domains, 
images, or both. Perhaps a good starting place to branch out would be the 
Knaster Buckethandle continuum (description in \cite{tree}). This 
indecomposable continuum can be neither the domain \cite{mio} nor the image 
\cite{jdm1} of a 2-to-1 (continuous) map, but in \cite{tree} an
example is given of an exactly 2-to-1 map from a hereditarily decomposable
tree-like continuum onto the Knaster Buckethandle continuum with exactly one 
discontinuity. However, the following is not known:
\vspace{.2in}
{\bf Question 11.} Is there a 2-to-1 finitely discontinuous function
defined on the Knaster Buckethandle continuum ?
\vspace{.2in}
There are a number of important facts, true about continuous k-to-1
functions, that remain true if finitely many discontinuities are allowed, even 
if the image is not required to be a continuum. For instance, (i)  the
dimension of the image is the dimension of the domain 
 for  continuous k-to-1 functions \cite{harrold2}, and, if the image
is compact,  for finitely discontinuous k-to-1 functions
\cite{fddendrite}; (ii) Gottschalk's result \cite{gotts} that no dendrite can 
be
the continuous k-to-1 image of any continuum is still true if the function
is allowed to have finitely many discontinuities \cite{fddendrite}; and
(iii) Harrold's original theorem in \cite{harrold1} that there is no continuous 

2-
to-1  map defined on $[0,1]$ also extends to the finitely discontinuous case,
\cite{reals}. But, oddly enough, the similar result, that the $n$-ball does not
support a continuous 2-to-1 map (Roberts \cite{roberts} for $ n = 2 $,
Civin \cite{civin} for  $n = 3 $   and \v{C}ernavskii \cite{cern2} for $ n > 2 
$),  is not true for finitely discontinuous functions. Krystyna 
Kuperberg constructed a 2-to-1 function defined on the unit square
that has exactly one discontinuity, and her example can be modified to make a 
2-to-1 function defined on the $n$-ball, for any $n > 2 $, with only
one discontinuity. This example has not appeared in print, and we will describe 

it here:
\vspace{.2in}
{\bf Example} (K. Kuperberg) A 2-to-1 function defined on the unit
disk with exactly one discontinuity.
\vspace{.2in}
We will use the following notation:
\begin{enumerate}
\item $D = \{(x,y)|x^{2} + y^{2} = 1\}$, the unit disk,
\item for each integer $n > 0$, $ S_{n} = \{(x,y)| (x + (n-1)/n)^{2} +
y^{2} = 1/n^{2} ;  y \geq 0 \} $, and
\item $S = D \setminus \cup^{\infty}_{n=1} S_{n} $.
\end{enumerate}
The domain of the function is some unit square. We will first remove a point 
$p$ from the boundary of this unit square and let $h$ denote a homeomorphism 
from  the square minus $p$ onto $S$. (Note: we are not suggesting that the 
boundary of $S$ is homeomorphic to the boundary of the square minus $p$.) 
We will, in the next paragraph, construct a continuous map $f$ on $S$ that is 
2-to-1 everywhere except for the one point $(0,-1)$ at the bottom of
$S$ at which it is 1-to-1. We will then extend the composition $ f \circ h $ to 

all of the square by 
mapping $p$ to $f((0,-1))$, to complete the construction of the function of 
the example. 
\vskip .2in
\centerline{\map}
\vskip .2in
The function $f$ on $S$ will be described as a series of identifications. 
First, for each point $(x,y)$ in the top half of $S$ (meaning $y > 0 $), 
identify $(x,y)$ and $(x,-y)$. Now the set of points of $S$ that have not been 
identified is the union of a countable collection $ \cal I $ of disjoint open 
intervals such that (1) each interval lies either on the $x$-axis or in the 
bottom half of $S$, (2) the sequence of intervals converges to the point $(-
1,0)$, and (3) the endpoints of each interval  are not in $S$. Next, locate  
the interval in $ \cal I $ containing $(0,-1)$ and identify 
each point $(x,y)$ in this interval with $(-x,y)$; thus the point $(0,-1)$ 
itself is not identified with another point. Then, for the intervals 
remaining in $ \cal I $, identify the first of these intervals with the second, 

the third of these intervals with the fourth, etc. We have now constructed $f$.

Note that although the original Kuperberg example does not have compact image, 
the image, $I$, can be made compact in the following way. First embed $I$ in 
its one point compactification and then identify the new point added with any 
point of $I$. The composition of these two maps is a one-to-one continuous 
function from $I$ onto a (compact) continuum and the composition of this 
composite function with Kuperberg's 2-to-1 function is again a 2-to-1 function 
with one discontinuity, but this time the image is a continuum.
 \section{Covering Maps.}
Covering maps defined on compact spaces are the  tamest of all k-to-
1  maps. Two-fold covering maps are related to {\em crisp} maps, i.e. maps
that are not just point-wise 2-to-1 but are continuum-wise 2-
to-1  in that, if $ C $ is a continuum in the image, then the inverse of $ C
$ consists of two disjoint continua each of which is mapped homeomorphically 
onto C. Every crisp map is a two-fold covering map and every two-fold covering 
map has a crisp restriction to a subcontinuum, \cite{hiimage}. So far as I 
know, this relationship has not been studied for integers greater than two, so 
a natural question is:
\vspace{.2in}
{\bf Question 12.} Define a map to be {\em $ k $-crisp } if for each continuum 
$ C $ in the image, the inverse of $ C $ consists of $ k $  disjoint continua, 
each of which is mapped homeomorphically onto C. What is the relationship, if 
any, between $ k $-crisp maps and $ k $-fold covering maps?
\vspace{.2in}

 \section{Simple Maps.}
In \cite{bors} K. Borsuk and R. Molski define {\em simple} maps to be 
continuous functions whose point inverses all have exactly one or two points. 
Simple maps share some of the strength of exactly 2-to-1 maps and
are much easier to construct. In fact, one instantly sees that the only space 
that does not support a simple map, that is not a homeomorphism, is the one 
point space. In case the simple map $ f $, defined on a compactum, is open, 
J. W. Jaworowski \cite{jaw} has shown that  $ f $ is  equivalent to a 
homeomorphism on the domain of period two. That is, the natural involution $ i 
$ on the compact domain defined by $ i(x) = x $ if $ f^{-1}f(x) $ has only one 
point and $ i(x) $ is the other point of $ f^{-1}f(x) $ otherwise, is a 
homeomorphism iff the simple map $ f $ is open. In contrast, if $ f $ is an 
open (exactly) 2-to-1 map, then $ f $ itself is locally one-to-one
and is a local homeomorphism \cite{hiimage}; but this is not true of simple 
maps (a simple example of this is folding an arc in half).

Exactly k-to-1 maps never change the dimension \cite{harrold2}, and
Jaworowski \cite{jaw} showed that open simple maps do not alter 
dimension; but in \cite{bors} Borsuk and Molski note that there is a simple map 

from the Cantor discontinuum onto an interval, so simple maps can raise 
dimension by one (but they point out that simple maps never change dimension 
other than to raise it by one). In a similar way, there is an natural simple 
map from the the Sierpi\'{n}ski universal curve into the plane that raises its 
dimension by one; however  W. D\c{e}bski and J. Mioduszewski have proved 
\cite{debskim1} the surprising result that every simple map from  the 
Sierpi\'{n}ski triangle into the plane has an image with empty interior (and so 

the image has dimension one). See \cite{debskim2} and \cite{debskim3} for other 

related results.

Borsuk and Molski \cite{bors} proved that every locally one-to-one map defined 
on a  compactum is a finite composition of simple maps. So in this sense, 
simple maps are building blocks for locally one-to-one maps on compacta. In 
\cite{siek} Sieklucki showed even more: Every map of finite order defined on a 
finite dimensional compactum  is a finite composition of simple maps. He also 
constructs an infinite dimensional counterexample. Since any finite composition 

of simple maps is necessarily of finite order, his theorem is the best possible 

for compacta. In response to the natural question of whether (or not) {\em 
open} maps of finite order (defined on a finite dimensional compacta) are 
finite compositions of {\em open} simple maps,  John Baildon \cite{bail} proved
that if $ f $ is an open simple map between 2-manifolds without boundaries, and 

if $ f $ is the composition of $ n $ open simple maps, then $ f $ has order $ 
2^{n} $. Hence, no such finite composition is possible for $ w = z^{3},  $ 
defined on the unit sphere, for instance. Note that Baildon adds to the 
definition of a simple map that it not be one-to-one.

These results do not extend, as is, to the exactly k-to-1 case. For
instance, there is a 3-to-1   map defined the unit interval onto a
simple closed curve, but it cannot be written as a composition of 2-to-
1  maps and 1-to-1  maps because there is no 2-to-1 map
defined on the unit interval at all, and finite compositions of one-to-one maps 

are homeomorphisms. But is there any kind of building block theory here?
\vspace{.2in}
{\bf Question 13.} Under what circumstances are k-to-1 maps finite
compositions of maps of lesser order?
\vspace{.2in}
 
\section{Definitions.}
\vskip .1in\noindent {\bf Adjacency matrix.} If $ V $ is a vertex set for a 
graph $H$, the
{\em adjacency matrix} is a matrix indexed by $ V  \times  V $ whose $ (v_{1}, 
v_{2}) $ entry is defined to be the number of edges in $H$ between $v_{1}$ and 
$v_{2}$. 
\vskip .1in\noindent {\bf Arc-continuum.} A continuum is an {\em arc-continuum} 

if each 
subcontinuum is either the whole continuum, a point or an arc.
\vskip .1in\noindent {\bf Arc-like.}  A continuum is {\em arc-like} if for each 

positive 
number $\epsilon$ there is an $\epsilon$-map from the continuum onto an arc, 
i.e. a  continuous function from the continuum onto an arc such each point 
inverse has diameter less than $\epsilon$.
\vskip .1in\noindent {\bf Confluent.} A function is {\em confluent} if for each 

continuum $C$ 
in the image, each component of the preimage of $C$ maps onto $C$.
\vskip .1in\noindent {\bf Continuum.}      A topological space is a {\em 
continuum} if it is 
connected, compact, and metric.
\vskip .1in\noindent {\bf Covering Map.} A continuous function $f$ from a space 

$ X$ onto a 
space $ Y $ is a {\em covering map} if for each point $y$ in $ Y $ there is an 
open set $ U $ containing $ y $ such that $ f^{-1}(U) $ is the union of 
finitely many disjoint open sets, each of which is mapped homeomorphically by 
$f$ onto $U$.
 \vskip .1in\noindent {\bf Crisp.} A map  $f$  is {\em crisp} if, for any 
proper subcontinuum 
$C$ of the image,  the inverse of $C$ is the union of two disjoint continua,  
each of which is mapped homeomorphically by $f$ onto $C$.
\vskip .1in\noindent {\bf Cut point.} A point $ x $ in a continuum $ X $ is a 
{\em cut point} 
if $ X \setminus \{x\} $ is not connected.
\vskip .1in\noindent {\bf Decomposable.} A continuum is {\em decomposable} if 
it is the union 
of two proper subcontinua.
\vskip .1in\noindent {\bf Dendrite.} A continuum is a {\em dendrite} if it is 
locally 
connected and contains no simple closed curve.
\vskip .1in\noindent {\bf Euler number.} The {\em Euler number} of a graph is 
the number of 
edges minus the number of vertices.
\vskip .1in\noindent {\bf Finitely discontinuous.} A function is {\em finitely 
discontinuous} 
if it has at most a finite number of discontinuities.
\vskip .1in\noindent {\bf Finite order.} A function has {\em finite order} if 
there is an 
integer $ k $ such that each point in the image has a preimage with no more 
than $ k $ points.
\vskip .1in\noindent {\bf Graph.} A continuum is a {\em graph} if it is 
homeomorphic to  the 
finite union of straight arcs and points.
\vskip .1in\noindent {\bf Hereditarily indecomposable.} A continuum is {\em 
hereditarily 
 indecomposable} if each subcontinuum is indecomposable.
\vskip .1in\noindent {\bf Indecomposable.} A continuum is {\em indecomposable}  

if it is not 
the  union of two proper (unequal to the continuum) subcontinua.
\vskip .1in\noindent {\bf Inverse Adjacency Matrix} If $f$ is a function from 
the vertex set of a 
graph $G$ onto the vertex set $V$ of a graph $H$, then the {\em inverse 
adjacency matrix} is indexed by $ V \times V$ and  its $ (v_{1}, v_{2}) $ entry 

is the number of edges in $G$ that go from any 
point of $f^{-1}(v_{1})$ to any point of $f^{-1}(v_{2})$.
\vskip .1in\noindent {\bf Involution.}  An {\em involution} is a function from 
a space into 
itself. It may or may not be continuous.
\vskip .1in\noindent {\bf k-to-1 } A function is {\em k-to-1 } if the
preimage of each point in the image has exactly $k$ points.
\vskip .1in\noindent {\bf k-crisp.} A map is {\em  k-crisp} if for each 
continuum $ C $
in the image, the inverse of $ C $ consists of $ k $ disjoint continua each of 
which is mapped homeomorphically onto $ C $.
\vskip .1in\noindent {\bf Local Cantor Bundle} A continuum is a {\em local 
Cantor bundle} if each 
point has a neighborhood homeomorphic to $ C \times (0,1) $, where $C$ denotes 
the Cantor discontinuum.
\vskip .1in\noindent {\bf Local Homeomorphism.} A function $f$ is a { \em local 

homeomorphism 
} if for each point $p$ in the domain, there is an open set $U$ containing $p$ 
such that  $f$  is a homeomorphism on $U$ and $f(U)$ is open.
\vskip .1in\noindent {\bf Map.} A function is a {\em map} if it is continuous.
\vskip .1in\noindent {\bf Non-orientable Arc-continuum.}  See "Orientable Arc-
continuum".
\vskip .1in\noindent {\bf Non-unicoherent.} A continuum is {\em non-
unicoherent} if it is  the 
union of two subcontinua whose intersection fails to be connected.
\vskip .1in\noindent {\bf Orientable Arc-continuum.} A general definition  can 
be found in 
\cite{over}, but for arc-continua that are local Cantor bundles, the definition 

is equivalent to the following natural one. The arc-continuum is {\em 
orientable} if each separate arc component can be parameterized (given a 
direction) so that no sequence of arcs going one direction converges to an arc 
going the other direction.
\vskip .1in\noindent {\bf Proper subcontinuum.} A subcontinuum of a continuum $ 
C $ is {\em proper} if it is not equal to $ C $. 
\vskip .1in\noindent {\bf Simple Map.} A continuous function is {\em simple} if 

each of its 
point inverses has cardinality 1 or 2.
\vskip .1in\noindent {\bf Tree.} A graph is a {\em tree} if it is connected and 

contains no 
simple closed curves.
 \vskip .1in\noindent {\bf Tree-like.}  A continuum is {\em tree-like} if for 
each positive 
number $\epsilon$ there is an  $\epsilon$-map from the continuum onto a tree. 
(See ``arc-like" for the definition of an  $\epsilon$-map.) 
\vskip .1in\noindent {\bf 2-to-1 } A function is {\em 2-to-1 } if the
preimage of each point in the image has exactly two points.
\vskip .1in\noindent {\bf Unicoherent.} See ``Non-unicoherent".
\vskip .1in\noindent {\bf Weak Confluence.} A function is {\em weakly confluent 

} if for each 
continuum $C$ in the image, at least one component of the preimage of $C$ maps 
onto $C$. 

\end{document}